\documentclass{article}
\usepackage{xurl}
\usepackage{tikz}
\usetikzlibrary{automata,positioning,arrows,fit}
\tikzstyle{morphism}=[-stealth, thick]
\tikzstyle{inclusion}=[{right hook}-stealth, thick]

\usepackage{amssymb}
\usepackage{amsmath}
\usepackage{phonetic}
\usepackage{xcolor}

\newcommand{\mName}[1]{\mathsf{#1}}
\newcommand{\mSet}[1]{\{ #1 \}}

\newcommand{\mAbs}[1]{\lvert #1 \rvert}

\newcommand{\mDot}{\mathrel{.}}

\newcommand{\cBoolTy}{\omicron}
\newcommand{\cB}{\cBoolTy}

\newcommand{\cFunTy}[2]{({#1} \rightarrow {#2})}
\newcommand{\cFunTyX}[2]{{#1} \rightarrow {#2}}

\newcommand{\cFunTyBX}[3]{\cFunTyX {#1} {\cFunTyX {#2} {#3}}}

\newcommand{\cFunTyCX}[4]{\cFunTyX {#1} {\cFunTyX {#2} {\cFunTyX {#3} {#4}}}}
\newcommand{\cProdTy}[2]{({#1} \times {#2})}

\newcommand{\cEqX}[2]{{#1} = {#2}}
\newcommand{\cFunApp}[2]{(#1\,#2)}
\newcommand{\cFunAppX}[2]{#1\,#2}

\newcommand{\cFunAppBX}[3]{\cFunAppX {\cFunAppX {#1} {#2}} {#3}}

\newcommand{\cFunAppCX}[4]{\cFunAppX {\cFunAppX {\cFunAppX {#1}{#2}}{#3}}{#4}}
\newcommand{\cFunAbs}[3]{(\lambda\, #1 : #2 \mDot #3)}
\newcommand{\cFunAbsX}[3]{\lambda\, #1 : #2 \mDot #3}

\newcommand{\cDefDesX}[3]{\mathrm{I}\, #1 : #2 \mDot #3}

\newcommand{\cAnd}[2]{({#1} \wedge {#2})}
\newcommand{\cAndX}[2]{{#1} \wedge {#2}}

\newcommand{\cImpliesX}[2]{{#1} \Rightarrow {#2}}

\newcommand{\cBin}[3]{({#1} \mathrel{#2} {#3})}
\newcommand{\cBinX}[3]{{#1} \mathrel{#2} {#3}}

\newcommand{\cBinBX}[5]{{#1} \mathrel{#2} {#3} \mathrel{#4} {#5}}

\newcommand{\cForall}[3]{(\forall\, #1 : #2 \mDot #3)}
\newcommand{\cForallX}[3]{\forall\, #1 : #2 \mDot #3}

\newcommand{\cForallBX}[5]{\forall\, #1 : #2,\, #3 : #4 \mDot #5}

\newcommand{\cForsomeX}[3]{\exists\, #1 : #2 \mDot #3}

\newcommand{\cSetTy}[1]{\mSet{#1}}

\newcommand{\cFracX}[2]{\frac {#1} {#2}}
\newcommand{\cAbs}[1]{\mAbs {#1}}

\newcommand{\cSum}[4]{\Big( \sum\limits_{{#1} = {#2}}^{#3} {#4} \Big)}
\newcommand{\cSumX}[4]{\sum\limits_{{#1} = {#2}}^{#3} {#4}}

\newcommand{\cLim}[3]{\Big( \lim\limits_{{#1} \to {#2}} {#3} \Big)}
\newcommand{\cLimX}[3]{\lim\limits_{{#1} \to {#2}} {#3}}

\newcommand{\cLimSeq}[2]{\Big( \lim\limits_{{#1} \to \infty} {#2} \Big)}
\newcommand{\cLimSeqX}[2]{\lim\limits_{{#1} \to \infty} {#2}}
\newcommand{\cIntegral}[4]{\Big( \int_{#1}^{#2} {#3}\,d{#4} \Big)}
\newcommand{\cIntegralX}[4]{\int_{#1}^{#2} {#3}\,d{#4}}

\newtheorem{thm}{Theorem}[section]

\newtheorem{thydef}[thm]{Theory Definition}
\newtheorem{thyext}[thm]{Theory Extension}
\newtheorem{indtypethyext}[thm]{Inductive Type Theory Extension}
\newtheorem{devdef}[thm]{Development Definition}
\newtheorem{devext}[thm]{Development Extension}
\newtheorem{thytransdef}[thm]{Theory Translation Definition}
\newtheorem{thytransext}[thm]{Theory Translation Extension}
\newtheorem{devtransdef}[thm]{Development Translation Definition}
\newtheorem{devtransext}[thm]{Development Translation Extension}
\newtheorem{deftransport}[thm]{Definition Transportation}
\newtheorem{thmtransport}[thm]{Theorem Transportation}
\newtheorem{grouptransport}[thm]{Group Transportation}

\newenvironment{theory-def}[5]
{
\color{brown!90!black}
\begin{thydef}[#1]\em
\noindent
\begin{itemize} \setlength{\itemsep}{0pt}
\item[]\hspace{-3ex}\textbf{Name:} #2
\item[]\hspace{-3ex}\textbf{Base types:} #3
\item[]\hspace{-3ex}\textbf{Constants:} #4
\item[]\hspace{-3ex}\textbf{Axioms:}
\end{itemize}
 #5
\end{thydef} 
}

\newenvironment{theory-ext}[6]
{
\color{brown!90!black}
\begin{thyext}[#1]\em
\noindent
\begin{itemize} \setlength{\itemsep}{0pt}
\item[]\hspace{-3ex}\textbf{Name:} #2
\item[]\hspace{-3ex}\textbf{Extends\ } #3
\item[]\hspace{-3ex}\textbf{New base types:} #4
\item[]\hspace{-3ex}\textbf{New constants:} #5
\item[]\hspace{-3ex}\textbf{New axioms:}
\end{itemize}
#6
\end{thyext}
}

\newenvironment{ind-type-theory-ext}[5]
{
\color{brown!90!black}
\begin{indtypethyext}[#1]\em
\noindent
\begin{itemize} \setlength{\itemsep}{0pt}
\item[]\hspace{-3ex}\textbf{Name:} #2
\item[]\hspace{-3ex}\textbf{Extends\ } #3
\item[]\hspace{-3ex}\textbf{New base type:} #4
\item[]\hspace{-3ex}\textbf{Constructors:}
\end{itemize}
#5
\end{indtypethyext}
}

\newenvironment{dev-def}[4]
{
\color{brown!90!black}
\begin{devdef}[#1]\em
\noindent
\begin{itemize} \setlength{\itemsep}{0pt}
\item[]\hspace{-3ex}\textbf{Name:} #2
\item[]\hspace{-3ex}\textbf{Bottom theory:} #3
\item[]\hspace{-3ex}\textbf{Definitions and theorems:}
\end{itemize}
#4
\end{devdef}
}

\newenvironment{dev-ext}[4]
{
\color{brown!90!black}
\begin{devext}[#1]\em
\noindent
\begin{itemize} \setlength{\itemsep}{0pt}
\item[]\hspace{-3ex}\textbf{Name:} #2
\item[]\hspace{-3ex}\textbf{Extends\ } #3
\item[]\hspace{-3ex}\textbf{New definitions and theorems:}
\end{itemize}
#4
\end{devext}
}

\newenvironment{theory-trans-def}[6]
{
\color{brown!90!black}
\begin{thytransdef}[#1]\em
\noindent
\begin{itemize} \setlength{\itemsep}{0pt}
\item[]\hspace{-3ex}\textbf{Name:} #2
\item[]\hspace{-3ex}\textbf{Source theory:} #3
\item[]\hspace{-3ex}\textbf{Target theory:} #4
\item[]\hspace{-3ex}\textbf{Base type mapping:}
\end{itemize}
#5
\begin{itemize}
\item[]\hspace{-3ex}\textbf{Constant mapping:}
\end{itemize}
#6
\end{thytransdef}
}

\newenvironment{theory-trans-ext}[7]
{
\color{brown!90!black}
\begin{thytransext}[#1]\em
\noindent
\begin{itemize} \setlength{\itemsep}{0pt}
\item[]\hspace{-3ex}\textbf{Name:} #2
\item[]\hspace{-3ex}\textbf{Extends\ } #3
\item[]\hspace{-3ex}\textbf{New source theory:} #4
\item[]\hspace{-3ex}\textbf{New target theory:} #5
\item[]\hspace{-3ex}\textbf{New base type mapping:}
\end{itemize}
#6
\begin{itemize}
\item[]\hspace{-3ex}\textbf{New constant mapping:}
\end{itemize}
#7
\end{thytransext}
}

\newenvironment{dev-trans-def}[6]
{
\color{brown!90!black}
\begin{devtransdef}[#1]\em
\noindent
\begin{itemize} \setlength{\itemsep}{0pt}
\item[]\hspace{-3ex}\textbf{Name:} #2
\item[]\hspace{-3ex}\textbf{Source development:} #3
\item[]\hspace{-3ex}\textbf{Target development:} #4
\item[]\hspace{-3ex}\textbf{Base type mapping:}
\end{itemize}
#5
\begin{itemize}
\item[]\hspace{-3ex}\textbf{Constant mapping:}
\end{itemize}
#6
\end{devtransdef}
}

\newenvironment{dev-trans-ext}[6]
{
\color{brown!90!black}
\begin{devtransext}[#1]\em
\noindent
\begin{itemize} \setlength{\itemsep}{0pt}
\item[]\hspace{-3ex}\textbf{Name:} #2
\item[]\hspace{-3ex}\textbf{Extends\ } #3
\item[]\hspace{-3ex}\textbf{New source development:} #4
\item[]\hspace{-3ex}\textbf{New target development:} #5
\item[]\hspace{-3ex}\textbf{New defined constant mapping:}
\end{itemize}
#6
\end{devtransext}
}

\newenvironment{def-transport}[9]
{
\color{brown!90!black}
\begin{deftransport}[#1]\em
\noindent
\begin{itemize} \setlength{\itemsep}{0pt}
\item[]\hspace{-3ex}\textbf{Name:} #2
\item[]\hspace{-3ex}\textbf{Source development:} #3
\item[]\hspace{-3ex}\textbf{Target development:} #4
\item[]\hspace{-3ex}\textbf{Development morphism:} #5
\item[]\hspace{-3ex}\textbf{Definition:}
\end{itemize}
#6
\begin{itemize}
\item[]\hspace{-3ex}\textbf{Transported definition:}
\end{itemize}
#7
\begin{itemize} \setlength{\itemsep}{0pt}
\item[]\hspace{-3ex}\textbf{New target development:} #8
\item[]\hspace{-3ex}\textbf{New development morphism:} #9
\end{itemize}
\end{deftransport}
}

\newenvironment{thm-transport}[8]
{
\color{brown!90!black}
\begin{thmtransport}[#1]\em
\noindent
\begin{itemize} \setlength{\itemsep}{0pt}
\item[]\hspace{-3ex}\textbf{Name:} #2
\item[]\hspace{-3ex}\textbf{Source development:} #3
\item[]\hspace{-3ex}\textbf{Target development:} #4
\item[]\hspace{-3ex}\textbf{Development morphism:} #5
\item[]\hspace{-3ex}\textbf{Theorem:}
\end{itemize}
#6
\begin{itemize}
\item[]\hspace{-3ex}\textbf{Transported theorem:}
\end{itemize}
#7
\begin{itemize}
\item[]\hspace{-3ex}\textbf{New target development:} #8
\end{itemize}
\end{thmtransport}
}

\newenvironment{group-transport}[9]
{
\color{brown!90!black}
\begin{grouptransport}[#1]\em
\noindent
\begin{itemize} \setlength{\itemsep}{0pt}
\item[]\hspace{-3ex}\textbf{Name:} #2
\item[]\hspace{-3ex}\textbf{Source development:} #3
\item[]\hspace{-3ex}\textbf{Target development:} #4
\item[]\hspace{-3ex}\textbf{Development morphism:} #5
\item[]\hspace{-3ex}\textbf{Definitions and theorems:}
\end{itemize}
#6
\begin{itemize}
\item[]\hspace{-3ex}\textbf{Transported definitions and theorems:}
\end{itemize}
#7
\begin{itemize} \setlength{\itemsep}{0pt}
\item[]\hspace{-3ex}\textbf{New target development:} #8
\item[]\hspace{-3ex}\textbf{New development morphism:} #9
\end{itemize}
\end{grouptransport}
}

\newcommand{\be}{\begin{enumerate}}
\newcommand{\ee}{\end{enumerate}}
\newcommand{\bi}{\begin{itemize}}
\newcommand{\ei}{\end{itemize}}
\newcommand{\bc}{\begin{center}}
\newcommand{\ec}{\end{center}}
\newcommand{\bsp}{\begin{sloppypar}}
\newcommand{\esp}{\end{sloppypar}}
\newcommand{\qzero}{${\cal Q}_0$}

\title{\bf An Alternative Approach to Formal Mathematics that Focuses on
  Communication and Accessibility}

\author{William M. Farmer\thanks{The author is a Professor in the
    Department of Computing and Software at McMaster University.  His
    email address is \texttt{wmfarmer@mcmaster.ca}.}  }

\date{July 7, 2026}

\begin{document}

\maketitle

\begin{abstract}

\emph{Formal mathematics} is mathematics done within the framework of
a formal logic.  It offers major benefits to mathematicians as well as
to computing professionals, engineers, and scientists who use
mathematics in their work.  The \emph{standard approach to formal
mathematics}, in which mathematics is done with the help of a proof
assistant and all details are formally proved and mechanically
checked, achieves these benefits and offers a very high level of
assurance that the results produced are correct.  However, since the
main goal of the standard approach is certification, the proof
assistants supporting the standard approach are generally complex,
based on unfamiliar logics, difficult to learn how to use, and far
removed from mathematical practice.  Thus the standard approach does
not adequately serve mathematics practitioners who are more interested
in communicating mathematical ideas than in formally certifying their
correctness or who prefer not to make the investment needed to gain
proficiency in the use of a proof assistant.

This paper presents an alternative to the standard approach that
focuses on communication and accessibility, the two weaknesses of the
standard approach.  It is called the \emph{free approach to formal
mathematics} since it is \emph{free} of the obligation to formally
prove and mechanically check all details of a mathematical
development.  The paper argues that the free approach would serve the
needs of the average mathematics practitioner much better than the
standard approach.  It describes an implementation of the free
approach based on a logic named \emph{Alonzo}, a practice-oriented
version of Alonzo Church's formulation of simple type theory.  And it
calls for the mathematics community to develop logics, software, and
libraries of formal mathematical knowledge to support the free
approach and to train mathematics practitioners to use them.

\end{abstract}

\section{Introduction}

\emph{Formal mathematics} is mathematics done within the framework of
a formal logic.\footnote{Formal mathematics is also called
``formalized mathematics'', but this latter term is often reserved for
the product that formal mathematics produces.}  Formal mathematics
differs from traditional mathematics in several ways.  In formal
mathematics all mathematical statements are expressed as formulas in a
formal language with a precise semantics, while in traditional
mathematics mathematical statements are expressed in a natural
language like English that does not have a precise semantics.  In
formal mathematics all assumptions and definitions are explicitly
stated, while in traditional mathematics many assumptions and
definitions are stated vaguely or not at all.  And in formal
mathematics all reasoning is performed in accordance with the logic's
notion of logical consequence, while in traditional mathematics
reasoning is often not grounded on a precise notion of logical
consequence.  Formal mathematics is based on ideas from logic, the
foundations of mathematics, and axiomatic systems.  These ideas were
developed over the last few centuries by mathematicians and logicians
such as (in chronological order) Leibniz, Lobachevsky, Boole,
Dedekind, Peirce, Cantor, Frege, Peano, Hilbert, Zermelo, Russell,
Curry, Tarski, Church, G\"{o}del, and Bourbaki.


Before the advent of electronic computers, interest in formal
mathematics was primarily a theoretical endeavor to develop a solid
foundation for mathematical thinking and to better understand what
mathematics is.  A salient example is the three-volume \emph{Principia
Mathematica}~\cite{WhiteheadRussell10} by Alfred North Whitehead and
Bertrand Russell that was intended to show how mathematical ideas
could be expressed in a formal logic and how the set-theoretic and
semantic paradoxes identified in the late 19th and early 20th
centuries could be resolved.

In the late 1960s, Nicolaas G. de Bruijn started the Automath
project~\cite{Nederpelt12} whose aim was the design of a formal
language in which mathematical ideas could be readily expressed and
checked by a computer.  The goal was to assist mathematicians in
developing new mathematics while remaining close to mathematical
practice.  The project was the first attempt to do formal mathematics
for practical purposes, and the Automath system was both the first
proof assistant\footnote{A \emph{proof assistant} is an interactive
software system for developing and checking formal proofs.} and the
first logical framework\footnote{A \emph{logical framework} is a
software system that provides the means to define and implement
multiple logics.}.  Following Automath, several other proof assistants
and logical frameworks have been developed including (in roughly
chronological order) Mizar~\cite{NaumowiczKornilowicz09}, Edinburgh
LCF~\cite{GordonEtAl79}, TPS~\cite{AndrewsEtAl96},
Nuprl~\cite{Constable86}, HOL~\cite{GordonMelham93},
Isabelle~\cite{Paulson94}, Rocq (formerly Coq)~\cite{Rocq},
ACL2~\cite{KaufmannEtAl00}, PVS~\cite{OwreEtAl96},
Metamath~\cite{Metamath}, IMPS~\cite{FarmerEtAl93,FarmerEtAl98b}, HOL
Light~\cite{Harrison09}, LEGO~\cite{Pollak94}, Agda~\cite{BoveEtAl09,
  Norell07}, Twelf~\cite{PfenningSchuermann99},
Lean~\cite{deMouraEtAl15}, and Naproche\cite{NaprocheWebSite}.

Although interest in using proof assistants to construct and
mechanically check formal proofs has been steadily
growing~\cite{BlanchetteMahboubi26}, formal mathematics is still not a
prominent component of mathematical practice.  We know of no studies
that have investigated the number of mathematics practitioners ---
mathematicians, computing professionals, engineers, and scientists ---
who work with formal mathematics.  However, a back-of-the-envelope
estimate is that no more (and probably far less) than 1\% of all
mathematics practitioners use a proof assistant in their
work.\footnote{There is very likely more than 30 million mathematics
practitioners in the world and many, many more if mathematics students
are included.  On the other hand, there is very likely fewer than
300,000 people in the world who use proof assistants.}  It is safe to
say that only a very small minority of mathematics practitioners are
pursuing formal mathematics today.

We believe formal mathematics holds great promise for mathematics.  In
particular, formal mathematics offers mathematics practitioners five
major benefits: (1) \emph{greater rigor}, (2) \emph{systematic
discovery of conceptual errors}, (3) \emph{software support},
(4)~\emph{mechanically checked results}, and (5) \emph{mathematical
knowledge as a formal structure} (see Section~\ref{sec:benefits} for
detailed descriptions of these benefits).  However, the promise of
formal mathematics has largely been unrealized since so few
mathematics practitioners ever do mathematics with the aid of a
formal~logic.

The purpose of this paper is to explain what formal mathematics is,
what benefits it offers, why it is unpopular, and finally and most
importantly, why an alternative approach is needed to make the
benefits of formal mathematics more widely available.

\section{What is Formal Mathematics?}

We have said that formal mathematics is mathematics done within the
framework of a formal logic.  So what does that mean?  To answer this
question we must first define a ``formal logic'', then define a
``theory'' of a formal logic, and finally define a ``theory morphism''
from one theory to another.

A \emph{formal logic} (\emph{logic} for short) is a family of
languages such that:

\be

  \item The languages of the logic have a \emph{precise common
  syntax}.  The languages have the same syntactic structure, but they
    have different vocabularies.  That is, they have the same logical
    symbols, but they have different nonlogical symbols.

  \item The languages of the logic have a \emph{precise common
  semantics with a notion of logical consequence}.  If $L$ is a
    language of the logic, $A$ is a sentence of~$L$, and $\Gamma$ is a
    set of sentences of $L$, then ``$A$ is a logical consequence of
    $\Gamma$'' means that in every possible situation, if all the
    sentences in $\Gamma$ are true, then $A$ is also true.

  \item There is a sound \emph{formal proof system} for the logic in
    which proofs can be syntactically constructed.  Let $L$ be a
    language of the logic, $A$ a sentence of $L$, and $\Gamma$ a set
    of sentences of $L$.  Soundness means that, if a proof of $A$ from
    $\Gamma$ can be constructed in the proof system, then $A$ is a
    logical consequence of $\Gamma$.

\ee 
This is a very general definition of a logic.  The key constituent of
a logic is its notion of logical consequence.  The definition covers
the various versions of first-order logic, simple type theory, and
dependent type theory as well as the various versions of set theory.
Note that a family of languages can be a family of one language by
itself.

Although the languages of the logic share a common syntax, the
expressions of a language can be presented, if desired, using multiple
notations including symbol-oriented notations like those found in
logic textbooks and controlled natural languages like the input
language for the Naproche proof assistant~\cite{NaprocheWebSite}.  The
proofs in a formal proof system can also be presented in multiple
notations including notations for describing formal proofs, notations
for prescribing how to construct formal proofs, and notations for
expressing formal proofs in controlled natural language.

A \emph{theory} of a logic is a pair $T = (L,\Gamma)$ where $L$ is a
language of the logic and $\Gamma$ is a set of sentences of~$L$.  The
members of $\Gamma$ are the \emph{axioms} of $T$; they serve as
background assumptions about the ideas that can be expressed in~$L$.
A theory is a fundamental unit of mathematical knowledge that is
\emph{developed} by defining new concepts in the theory and by posing
and then proving conjectures in the theory.  Mathematical knowledge
can be organized either as a single, large, highly developed theory or,
much better, as a network of interconnected theories of varying size
and development.  We will look next at the mechanism for connecting
theories.

Let $T$ and $T'$ be theories of a logic. A \emph{theory morphism from
$T$ to $T'$} is a structure-preserving translation $\Phi$ from the
expressions of $T$ to the expressions of $T'$ that is also meaning
preserving in the sense that each theorem of~$T$ is translated to a
theorem of $T'$.  The concepts (established by axioms and definitions)
and facts (established by theorems) of $T$ can be transported to~$T'$
via $\Phi$.  For example, a theory morphism $\Phi_{\mathsf{MON}
  \mapsto \mathsf{COF}^+}$ from a theory $\mathsf{MON}$ of monoids to
a theory $\mathsf{COF}$ of complete ordered fields (which specifies
real number mathematics) that maps the binary operator $\cdot$ and the
identity element $\mathsf{e}$ in $\mathsf{MON}$ to $+$ and $0$,
respectively, in $\mathsf{COF}$ can be used to transport the
theorem~\textsf{id-elt-is-unique} of $\mathsf{MON}$ that
says~$\mathsf{e}$ is the unique identity element with respect to
$\cdot$ to the theorem of $\mathsf{COF}$ that says $0$ is the unique
identity element with respect to~$+$.  That is, $\Phi_{\mathsf{MON}
  \mapsto \mathsf{COF}^+}$ translates a definition or theorem (which
can be an axiom) $X$ of~$\mathsf{MON}$ to an instance of $X$ that is a
definition or theorem of $\mathsf{COF}$.  (Note that the theory
morphism $\Phi_{\mathsf{MON} \mapsto \mathsf{COF}^*}$ from
$\mathsf{MON}$ to $\mathsf{COF}$ that translates the binary operator
$\cdot$ and the identity element $\mathsf{e}$ to $*$ and~$1$,
respectively, can transport \textsf{id-elt-is-unique} to the theorem
of $\mathsf{COF}$ that says~$1$ is the unique identity element with
respect to~$*$.)  Theory morphisms are thus unidirectional conduits
through which information in the form of concepts and facts flows from
usually more abstract to usually more concrete theories.

The \emph{little theories method}~\cite{FarmerEtAl92b} is a method for
organizing mathematical knowledge as a network of interconnected
theories called a \emph{theory graph}~\cite{KohlhaseEtAl10}.  The
nodes of the graph are \emph{theories} of a logic and the directed
edges are \emph{theory morphisms} between the theories.  To maximize
clarity, each mathematical topic is developed in the ``little theory''
in the theory graph that has the most convenient vocabulary and level
of abstraction.  To minimize redundancy, the definitions and theorems
produced in this little theory are transported, as needed, to other
theories via the theory morphisms in the theory graph instead of
directly producing them in the theories where they are~needed.

Figure~\ref{fig:monoid-theory} shows an example of a theory graph
associated with monoid theory that is taken
from~\cite{FarmerZvigelsky25} (which we will discuss in
Section~\ref{sec:example}).  A theory morphism that is an inclusion
(i.e., an identity mapping) is designated by a $\hookrightarrow$ arrow
and a noninclusion is designated by a $\rightarrow$ arrow.  Note that,
when there is an inclusion from $T$ to $T'$, then $T'$ is an
\emph{extension} of $T$ (i.e., the language of $T$ is a sublanguage of
the language of $T'$ and the axioms of $T$ are a subset of the axioms
of $T'$).  There are many, many more useful theory morphisms that are
not shown in Figure~\ref{fig:monoid-theory}, including
$\Phi_{\mathsf{MON} \mapsto \mathsf{COF}^+}$ (and $\Phi_{\mathsf{MON}
  \mapsto \mathsf{COF}^*}$) as well as a vast number of other theory
morphisms to the theory~\textsf{COF}.

\begin{figure*}[t]
\tikzset{every loop/.style={min distance=5mm,in=0,out=60,looseness=6}}
\tikzset{state/.style={circle, draw, minimum size=2cm}}
\tikzset{node distance=3.4cm,on grid,auto}
\center
\begin{tikzpicture}
  \node[state, thick] (mon) {\textsf{MON}};
  \node[state, thick, align=center] (com-mon) [above = of mon] 
    {\footnotesize{\textsf{COM-}} \\ \footnotesize{\textsf{MON}}};
  \node[state, thick, align=center] (com-hom) [left = of mon] 
    {\footnotesize{\textsf{MON-}} \\ \footnotesize{\textsf{HOM}}};
  \node[state, thick, align=center] (com-mon-cof) [right = of com-mon] 
    {\scriptsize{\textsf{COM-MON-}} \\ \scriptsize{\textsf{over-COF}}};
  \node[state, thick, align=center] (mon-cof) [right = of mon]  
    {\footnotesize{\textsf{MON-}} \\ \footnotesize{\textsf{over-COF}}};
  \node[state, thick, align=center] (mon-act) [below = of mon] 
    {\footnotesize{\textsf{MON-}} \\ \footnotesize{\textsf{ACT}}};
  \node[state, thick, align=center] (one) [left = of mon-act] 
    {\footnotesize{\textsf{ONE-}} \\ \footnotesize{\textsf{BT}}}; 
  \node[state, thick, align=center] (fun-comp) [below = of one]  
    {\footnotesize{\textsf{FUN-}} \\ \footnotesize{\textsf{COMP}}};
  \node[state, thick, align=center] (one-sc) [below = of mon-act]
    {\footnotesize{\textsf{ONE-BT-}} \\ \footnotesize{\textsf{with-SC}}};
  \node[state, thick] (cof) [right = of one-sc] {\textsf{COF}};
  \node[state, thick, align=center] (com-mon-act-cof) [right = of com-mon-cof] 
     {\scriptsize{\textsf{COM-}} \\ 
      \scriptsize{\textsf{MON-ACT-}} \\
      \scriptsize{\textsf{over-COF}}};
  \node[state, thick] (str) [right = of cof] {\textsf{STR}};
  \path [inclusion]
    (mon)             edge (com-mon)
                      edge (mon-cof)
                      edge (mon-act)
    (com-mon)         edge (com-mon-cof)
    (com-mon-cof)     edge (com-mon-act-cof)
    (mon-act)         edge [bend right] (com-mon-act-cof)
    (mon-cof)         edge (com-mon-cof)
    (one)             edge (one-sc)
    (cof)             edge (str)
                      edge (mon-cof);
  \path [morphism]
    (mon)             edge [out=115, in=135, looseness=15] (mon)
                      edge [out=105, in=145, looseness=10] (mon)
                      edge [out=170, in=30] (com-hom)
                      edge [out=190, in=10] (com-hom)
                      edge (one)
    (com-hom)         edge [bend right] (mon)
    (mon-act)         edge (one-sc)
                      edge [bend right] (mon)
    (fun-comp)        edge (one)
    (mon-cof)         edge [bend left] (str);
\end{tikzpicture}
\caption{A Theory Graph for Monoid Theory}
\label{fig:monoid-theory}
\end{figure*}

Revisiting our example above, the theorem \textsf{id-elt-is-unique}
about the uniqueness of the identity element in a monoid is stated and
proved in $\mathsf{MON}$ and then transported to $\mathsf{COF}$ via
$\Phi_{\mathsf{MON} \mapsto \mathsf{COF}^+}$ instead of directly stating
and proving the instance of this theorem in $\mathsf{COF}$ or in any
other theory that embodies a monoid structure.  The little theories
method thus provides a very strong form of \emph{polymorphism}: Any
definition or theorem $X$ produced in a theory $T$ can be reused in
every theory $T'$ for which there is a theory morphism from $T$ to
$T'$.

Now we can illustrate with a simple example how ``mathematics is done
within a formal logic''.  We will use the little theories method.
Suppose that we would like to show that a conjecture~$C$ is true
within some logic.  Let $G$ be a theory graph of the logic that has
been constructed using the little theories method.  We proceed by
obtaining a theory $T = (L,\Gamma)$ that best captures the context in
which we want to consider~$C$.  We obtain $T$ by either selecting a
theory that is already in~$G$ or by constructing a new theory that is
added to~$G$ and connected as needed to other theories in $G$ via
theory morphisms.  Then we express $C$ as a sentence $A_C$ of the
language $L$ augmented with the definitions made in $T$ and the
definitions transported to~$T$ from other theories via theory
morphisms.  And finally we try to prove $A_C$ is a theorem of $T$ by
showing $A_C$ is a logical consequence of the axioms of~$T$.  This can
be done either by composing a traditional-style proof that $A_C$
follows from $\Gamma$ (the axioms of~$T$), the theorems proved in $T$,
and the theorems transported to $T$ via theory morphisms or by
constructing a formal proof in the formal proof system of the logic.

We thus see that doing formal mathematics involves the following five
activities:

\be

  \item Constructing theories in a logic.

  \item Defining new concepts in a theory.

  \item Posing and then proving conjectures in a theory.

  \item Constructing theory morphisms.

  \item Transporting definitions and theorems via theory morphisms.

\ee
(The last two activities may not be needed if all the work is done in
one theory.)

\section{What are the Benefits of Formal Mathematics?}\label{sec:benefits}

There are five major benefits of formal mathematics.

The first benefit is that \emph{mathematics can be done with greater
rigor}.  A prime goal of mathematics is to express and reason about
mathematical ideas in a highly rigorous manner.  There are two aspects
of traditional mathematical practice that stand in the way of rigor.
The first is that mathematical ideas are usually expressed in a
natural language, such as English, that does not have a precise
semantics, and the second is that assumptions, definitions, theorems,
and even the notion of logical consequence are often expressed
imprecisely or implicitly.  These shortcomings do not arise in formal
mathematics since (1) a logic has a precise notion of logical
consequence, (2)~all mathematical concepts and statements are
expressed in a language that has a precise semantics, and (3)~all
assumptions, definitions, and theorems are expressed explicitly as
axioms, definitions, and theorems of a theory.

The second benefit is that \emph{conceptual errors can be
systematically discovered}.  It is easy to make conceptual errors in
mathematics, especially in mathematics that is complex or not well
understood.  A logic offers a precise conceptual framework to a
mathematics practitioner.  The process of expressing ideas in the
logic can systematically lead to the discovery of many conceptual
errors.  Let us consider the following very simple example.  Suppose
that a definition is employed early in a mathematical development, and
a closely related but inequivalent definition is employed by mistake
late in the development in place of the first definition.  This
mistake, if not discovered, could lead to invalid results.  If the
development were done within a formal logic, this mistake would be
easily caught since every definition would be expressed as a sentence
with a precise syntax and semantics.  The process of expressing
mathematical ideas in a formal logic naturally leads to many
conceptual errors being caught similarly to how type errors are caught
in a modern programming language by type checking.  Thus conceptual
errors can be discovered systematically in formal mathematics in a way
that is largely not possible in traditional mathematics.

The third benefit is that \emph{mathematics can be done with software
support}.  Since the languages of a formal logic have a precise common
syntax, the expressions of a language can be represented as data
structures.  The expressions can then be analyzed, manipulated, and
processed by software applied to the data structures that represent
them.  This, in turn, enables the study, discovery, communication, and
certification of mathematics to be done with the aid of software.  The
software can provide many useful services such as performing numerical
and symbolic computations, solving various kinds of problems, and
displaying mathematical data and expressions.  Since the languages
also have a precise common semantics, there is a precise basis for
verifying the correctness of this software.

The fourth benefit is that \emph{results can be mechanically checked}.
Formal proofs can be represented as data structures, and software can
be used to check that these data structures represent actual proofs in
the formal proof system of the logic.  \emph{Proof assistants} can be
used to help humans construct formal proofs, and \emph{automated
theorem provers} can be used to automatically discover formal proofs.
Since the software needed to check the correctness of the formal
proofs is often very simple and easily verified itself, it is possible
to verify the correctness of the formal proofs with a very high level
of assurance.  Thus mathematical results can be mechanically checked
by constructing formal proofs of them in the formal proof system of
the logic and then using software to check the correctness of the
formal proofs.  High assurance of mathematical correctness is
especially needed for cutting-edge mathematics
research~\cite{Avigad24} and the development of safety-, security-,
and mission-critical software~\cite{BlanchetteMahboubi26}.

Finally, the fifth benefit is that \emph{we can regard mathematical
knowledge as a formal structure}.  Using the little theories method, a
body of mathematical knowledge can be represented as a theory graph
built by creating theories, defining new concepts, posing and then
proving conjectures, constructing theory morphisms, and transporting
definitions and theorems via theory morphisms.  The knowledge embodied
in a formal structure of this kind can be studied, developed,
searched, and presented using software.

The benefits of formal mathematics are huge.  Greater rigor and
discovering conceptual errors have been principal goals of
mathematicians for thousands of years.  Software support can greatly
extend the reach and productivity of mathematics practitioners.
Mechanically checked results can drive mathematics forward in areas
where the ideas are poorly understood (often due to their novelty),
highly complex, or underlying a critical real-world application.  And
mathematical knowledge as a formal structure can enable the techniques
and tools of mathematics and computing to be applied to mathematical
knowledge~itself.

\section{Why is Formal Mathematics Unpopular?}

One would expect that, by virtue of these five benefits, formal
mathematics would be very popular among mathematics practitioners.
Yet the opposite is true.  So why is formal mathematics so unpopular?

The answer lies in the \emph{standard approach to formal mathematics}
in which mathematics is done with a proof assistant and all details
are formally proved and mechanically checked.  The standard approach
has three major strengths:

\be

  \item It achieves all five benefits of formal mathematics mentioned
    above.

  \item All theorems are verified by mechanically checked formal
    proofs.  Thus there is a very high level of assurance that the
    results produced are correct.

  \item \bsp There are several powerful proof assistants available,
    such as ACL2, HOL, HOL Light, Isabelle/HOL, Lean, Metamath/ZFC,
    Mizar, PVS, and Rocq (formerly Coq), that support the
    approach. \esp

\ee

\noindent
But the standard approach also has two important weaknesses:

\be

  \item It prioritizes certification over communication.  

  \item It is not accessible to the great majority of mathematics
    practitioners.  

\ee

Proof assistants are designed primarily for formally \emph{certifying
mathematical results} and only secondarily for \emph{communicating
mathematical ideas}.  Formal proofs are great for certifying that
something is true but often do a poor job of communicating why
something is true.  Most mathematics practitioners --- including
mathematicians --- are much more interested in communicating
mathematical ideas than in formally certifying their correctness.
This is particularly true for applications that involve
well-understood mathematics, the kind of mathematics that arises in
mathematics education and routine applications.

To apply the standard approach, a mathematics practitioner needs a
solid understanding of formal logic and how to use a proof assistant
to develop theories in the underlying logic of the proof assistant.
Since the focus of the standard approach is on certification, the
proof assistants supporting the approach are generally complex, based
on unfamiliar logics, difficult to learn how to use, and far removed
from mathematical practice.  The investment needed to learn a specific
proof assistant is so high that the user is usually tied to just one
proof assistant and is proficient with just one logic, one formal
proof system, one set of notation, and one set of software tools.
This makes it difficult for a proof assistant user to share their
results with users of other proof assistants.  Having to express and
check all details in an often unfamiliar logic using a complex proof
assistant that utilizes certification-oriented ways of expressing and
reasoning about mathematics is a bridge too far for many mathematics
practitioners.

Thus the standard approach to formal methods does not adequately serve
the average mathematics practitioner.  As a result, the average
mathematics practitioner does not have access to the benefits of
formal mathematics since the standard approach is pretty much the only
game in town.  

We are not saying that proof assistant developers care only about
certification.  They want to improve support in their systems for
communication and to make their systems more accessible.  Significant
progress has been made in both directions, and this has helped
accelerate the growth of the proof assistant user community.  There
has been particularly strong recent growth in the number of people
using the Lean proof assistant and working on the development of
Lean's Mathlib library of mathematical knowledge and programming
infrastructure.  Nevertheless, nearly all contemporary proof
assistants prioritize certification over communication and
accessibility.  (The Naproche system is one notable exception.)


\section{Why is an Alternative Approach Needed?}

We strongly believe an alternative to the standard approach to formal
mathematics is needed that focuses on two goals, \emph{communication}
and \emph{accessibility}, the two weaknesses of the standard approach.
We propose an alternative approach of this kind called the \emph{free
approach to formal mathematics} since it is \emph{free} of the
obligation to formally prove and mechanically check all details of a
mathematical development using a proof assistant.  To achieve the
goals of communication and accessibility, an implementation of the
free approach needs to satisfy the following requirements:
\be

  \item[R1.] \emph{The underlying logic is fully formal and supports
  standard mathematical practice.}  Supporting mathematical practice
    makes the logic easier to learn and use and makes formalization a
    more natural process.  It also makes formal mathematical statements
    easier to read, write, and reason about --- which helps the
    developer to identify mistakes and see connections.

  \item[R2.] \emph{Proofs can be traditional, formal, or a combination
  of the two.}  This flexibility in how proofs are written enables
    proofs to be a vehicle for communication as well as certification.
    Traditional proofs are usually easier to read and write and better
    suited for communicating the ideas behind proofs than pure formal
    proofs.

  \item[R3.] \emph{There is support for organizing mathematical
  knowledge using the little theories method.}  This enables
    mathematical knowledge to be formalized to maximize clarity and
    minimize redundancy (as explained above).

  \item[R4.] \emph{There are several levels of supporting software.}
    The levels can range from just LaTeX support to a full proof
    assistant to an \emph{interactive mathematics laboratory
    (IML)}~\cite{Farmer25} that provides support for inference as well
    as other aspects of mathematics~\cite{CaretteEtAl21}.  An
    important intermediate example is a software system that is the
    same as a proof assistant except it only supports the development
    of traditional proofs.  A system of this kind, a
    \emph{traditional-proof proof assistant}, would be much simpler
    and thus easier to both implement and learn how to use than a
    typical proof assistant.  The user can thus choose the level of
    software support they want to have and the level of investment in
    learning the software they want to make.

\ee
These four requirements characterize the free approach.

The free approach is intended to achieve all five benefits of
formal mathematics mentioned above, but it cannot achieve the same
level of assurance, as the standard approach, that the results
produced are correct.  This is because, in order to prioritize
communication and accessibility over certification, the free
approach does not require that all details are formally proved and
mechanically checked.

An implementation of the free approach --- with support for standard
mathematical practice, traditional proofs, the little theories method,
and several levels of software --- is likely to serve the needs of the
average mathematics practitioner much better than an implementation of
the standard approach.  This is especially true when the mathematical
knowledge involved is well understood and certification via
traditional proof is adequate for the purpose at~hand.

In summary, the free approach is not a replacement for the standard
approach, but we believe it would be more useful, accessible, and
natural than the standard approach for the great majority of
mathematics practitioners.  Thus it would bring the benefits of formal
mathematics to vastly more people than the standard approach.

\section{Other Alternative Approaches}

Other alternative approaches to formal mathematics have been proposed.
Three of the most notable reside in the space between traditional
mathematics and fully certified formal mathematics.

The first is Tom Hales' \emph{formal abstracts in mathematics}
project~\cite{Hales18} in which proof assistants are used to create
\emph{formal abstracts}, which are formal presentations of
mathematical theorems without formal proofs.  The formal abstracts
approach holds open the possibility of adding formal proofs later.
The free approach differs from the formal abstracts approach in that
it is not committed to using proof assistants and is not required to
facilitate the addition of formal proofs.

The second is Michael Kohlhase's \emph{flexiformal
mathematics}~\cite{Iancu17,Kohlhase12} initiative in which mathematics
is a mixture of traditional mathematics with traditional proofs and
formal mathematics with formal proofs.  The formal mathematics can
also be presented simultaneously in both traditional and formal forms.
A goal of the flexiformal mathematics approach is to gradually
transform the traditional mathematics into formal mathematics.  In
contrast to the flexiformal mathematics approach, the free approach
produces mathematics that is fully formal except that proofs may be
traditional.

The third is the \emph{natural language approach} in which mathematics
is performed using a proof assistant like
Naproche~\cite{NaprocheWebSite} or Natty~\cite{Dingle24} that accepts
a controlled form of mathematical English (or another natural
language) as input.  The natural language approach produces fully
certified formal mathematics written in natural language.  The use of
natural language can improve communication, but the natural language
approach is a really just an instance of the standard approach in
which natural language notation is used to present formal expressions.
It is thus, unlike the free approach, a certification-oriented
approach to formal mathematics.

\section{An Implementation of the Free Approach}

We have developed an implementation of the free approach based on
Alonzo~\cite{Farmer25}, a practice-oriented classical higher-order
predicate logic that extends first-order logic.  Named in honor of
Alonzo Church, Alonzo is a version of Church's type
theory~\cite{SEP-ChurchsTypeTheory25,Church40}, Church's formulation
of simple type theory~\cite{Farmer08}. It is closely related to Peter
Andrews' {\qzero}~\cite{Andrews02} and
LUTINS~\cite{Farmer90,Farmer94}, the logic of the IMPS proof
assistant~\cite{FarmerEtAl93}.

Alonzo is designed to be as close to mathematical practice as
possible.  By virtue of being a form of predicate logic, it is built
on ideas that are widely familiar to mathematics practitioners.  It is
a type theory, but its type system is simpler than the type systems of
most of the type theories employed in proof assistants.  It is
equipped with facilities, including higher-order quantification and
definite description, for reasoning about functions, sets, tuples,
lists, and mathematical structures.  Unlike traditional predicate
logics, Alonzo admits partial functions and undefined expressions in
accordance with the approach employed in mathematical practice that we
call the \emph{traditional approach to
undefinedness}~\cite{Farmer04,Farmer25}.  With this approach,
non-Boolean expressions that naturally have no value --- like the top
of an empty stack or $\mathsf{lim}_{x \to 0}\,\mathsf{sin}(1/x)$ ---
are considered \emph{undefined} and denote nothing at all, but Boolean
expressions --- even if they contain undefined non-Boolean expressions
--- always denote either true or false.\footnote{The traditional
approach to undefinedness is exemplified in Michael Spivak's famous
textbook \emph{Calculus}; see the discussion in~\cite{Farmer04}.}

Alonzo has a simple syntax with a \emph{formal notation} for machines
and a \emph{compact notation} for humans that closely resembles the
notation found in mathematical practice.  The compact notation is
defined by the extensive set of \emph{notational definitions and
conventions} given in~\cite{Farmer25}.  Alonzo has two semantics, one
for mathematics based on \emph{standard models} and one for logic
based on Henkin-style \emph{general models}.  Taken as a whole, it
offers a very good balance between \emph{theoretical expressivity}
(the measure of what ideas can be expressed without regard to how the
ideas are expressed) and \emph{practical expressivity} (the measure of
how readily ideas can be expressed).

Figure~\ref{fig:COF-calculus} shows part of the development of the
Alonzo theory $\mathsf{COF}$ of complete ordered fields taken from
Chapter 13 (Real Number Mathematics) of~\cite{Farmer25}.  It contains
several definitions and theorems from calculus presented using various
notational definitions including those given in
Figure~\ref{fig:notational-defs}.  $R$ is the type of real numbers and
$N_{\cSetTy {R}}$ denotes the set of natural numbers as a subset of
the set of real numbers.  The reader should note that the functions in
these definitions and theorems are ``Curryed''. For example, a binary
function that takes values of type $\alpha$ and~$\beta$ as input and
returns a value of type $\gamma$ as output is represented in Curryed
form as a function of type ${\cFunTyX {\alpha} {\cFunTy {\beta}
    {\gamma}}}$ or, more simply, as ${\cFunTyBX {\alpha} {\beta}
  {\gamma}}$.  (A binary function could also be represented in
Alonzo, if desired, as a function of type ${\cFunTyX {\cProdTy
    {\alpha} {\beta}} {\gamma}}$.)  The reader should also note that
the definitions of a limit of a function and a limit of a sequence
employ definite descriptions of the form ${\cDefDesX {b} {R}
  {\textbf{A}_\cB}}$ that denotes the unique real number $b$ that
satisfies the formula $\textbf{A}_\cB$ if such a real number exists
and is undefined otherwise.

The definitions and theorems presented in
Figure~\ref{fig:COF-calculus} utilize the strong support Alonzo
provides for partial functions and undefined expressions.  Partial and
total functions are handled in Alonzo in exactly the same way; they
cannot be distinguished in Alonzo syntactically.  For example, the
function application ${\cFunAppX {\mName{lim}_{\cFunTyBX {\cFunTy {R}
        {R}} {R} {R}}} {\mathbf{F}_{\cFunTy {R} {R}}}}$ (where the
constant ${\mName{lim}_{\cFunTyBX {\cFunTy {R} {R}} {R} {R}}}$ is
defined by $\mName{Def13}$) can denote either a partial or total
function, depending on what function ${\mathbf{F}_{\cFunTy {R} {R}}}$
denotes.  Since the application of equality and other predicates to
undefined expressions always denotes false, undefinedness conditions
can be expressed implicitly.  For example, in $\mName{Def14}$ the
equality ${\cEqX {\cLimX {x} {a} {\cFunAppX {f} {x}}} {\cFunAppX {f}
    {a}}}$ asserts (1) the limit ${\cLimX {x} {a} {\cFunAppX {f}
    {x}}}$ is defined, (2) the function application ${\cFunAppX {f}
  {a}}$ is defined, and (3) the values of ${\cLimX {x} {a} {\cFunAppX
    {f} {x}}}$ and ${\cFunAppX {f} {a}}$ are equal.  Thus statements
that include definedness conditions can usually be expressed in Alonzo
very concisely.

We believe that, by virtue of its syntax, semantics, and notational
definitions and conventions, Alonzo satisfies requirement R1 as well
or better than almost any other formal logic.

\begin{figure*}
\bc
\begin{tabular}{|l|}
\hline

  \phantom{x}\\

  $\mName{Def10}$: ${\cEqX {\mName{lim}_{\cFunTyBX {\cFunTy {R} {R}} {R} {R}}} {}}$\\
    \hspace*{2ex}${\cFunAbsX {f} {\cFunTyX {R} {R}} {\cFunAbsX {a} {R} 
    {\cDefDesX {b} {R} {}}}}$\\
    \hspace*{4ex}$({\cForallX {e} {R} {\cImpliesX {\cBinX {0} {<} {e}} {}}}$\\
    \hspace*{6ex}$({\cForsomeX {d} {R} {\cAndX {\cBinX {0} {<} {d}} {}}}$\\
    \hspace*{8ex}$({\cForallX {x} {R} {\cImpliesX
    {\cBinBX {0} {<} {\cAbs {\cBinX {x} {-} {a}}} {<} {d}} {}}}$\\
    \hspace*{10ex}${\cBinX {\cAbs {\cBinX {\cFunAppX {f} {x}} {-} {b}}} {<} {e}})))$
    \hfill (limit of a function).\\[2ex]

  $\mName{Def13}$: ${\cEqX {\textsf{lim-seq}_{\cFunTyX {\cFunTy {R} {R}} {R}}} {}}$\\
    \hspace*{2ex}${\cFunAbsX {s} {\cFunTyX {N_{\cSetTy {R}}} {R}} 
    {\cDefDesX {b} {R} {}}}$\\
    \hspace*{4ex}$({\cForallX {e} {R} {\cImpliesX {\cBinX {0} {<} {e}} {}}}$\\
    \hspace*{6ex}$({\cForsomeX {m} {N_{\cSetTy {R}}} {}}$\\
    \hspace*{8ex}$({\cForallX {n} {N_{\cSetTy {R}}} {\cImpliesX {\cBinX {m} {<} {n}} {}}}$\\
    \hspace*{10ex}${\cBinX {\cAbs {\cBinX {\cFunAppX {s} {n}} {-} {b}}} {<} {e}})))$
    \hfill (limit of a sequence).\\[2ex]

  $\mName{Def14}$: ${\cEqX {\textsf{cont-at}_{\cFunTyBX {\cFunTy {R} {R}} {R} {\cB}}} {}}$\\
    \hspace*{2ex}${\cFunAbsX {f} {\cFunTyX {R} {R}} {\cFunAbsX {a} {R}}
    {\cEqX {\cLimX {x} {a} {\cFunAppX {f} {x}}} {\cFunAppX {f} {a}}}}$ 
    \hfill (continuous at a point).\\[2ex]

  $\mName{Def18}$: 
    ${\cEqX {\textsf{cont-on-closed-int}_{\cFunTyCX {\cFunTy {R} {R}} {R} {R} {\cB}}} {}}$\\
    \hspace*{2ex}${\cFunAbsX {f} {\cFunTyX {R} {R}} 
    {\cFunAbsX {a} {R}} {\cFunAbsX {b} {R}} {}}$\\
    \hspace*{4ex}${\cAndX
    {\cForall {x} {R} {\cImpliesX 
    {\cBinBX {a} {<} {x} {<} {b}}
    {\cFunAppBX {\textsf{cont-at}} {f} {x}}}} {}}$\\
    \hspace*{4ex}${\cAndX {\cFunAppBX {\textsf{right-cont-at}} {f} {a}} {}}$\\
    \hspace*{4ex}${\cFunAppBX {\textsf{left-cont-at}} {f} {b}}$
    \hfill (continuous on a closed interval).\\[2ex]

  $\mName{Def19}$: ${\cEqX 
    {\textsf{deriv-at}_{\cFunTyBX {\cFunTy {R} {R}} {R} {R}}} {}}$\\
    \hspace*{2ex}${\cFunAbsX {f} {\cFunTyX {R} {R}} {\cFunAbsX {a} {R}
    {\cLimX {h} {0} {\cFracX {\cBinX {\cFunAppX {f} {\cBin {a} {+} {h}}} {-} 
    {\cFunAppX {f} {a}}} {h}}}}}$ \hfill (derivative at a point).\\[2ex]

  $\mName{Def22}$: ${\cEqX 
    {\mName{deriv}_{\cFunTyX {\cFunTy {R} {R}} {\cFunTy {R} {R}}}}
    {\cFunAbsX {f} {\cFunTyX {R} {R}} {\cFunAbsX {x} {R}
    {\cFunAppBX {\textsf{deriv-at}} {f} {x}}}}}$
    \hfill (derivative).\\[2ex]

  $\mName{Def26}$: ${\cEqX 
    {\mName{integral}_{\cFunTyCX {\cFunTy {R} {R}} {R} {R} {R}}} {}}$\\
    \hspace*{2ex}${\cFunAbsX {f} {\cFunTyX {R} {R}} 
    {\cFunAbsX {a} {R}} {\cFunAbsX {b} {R}} {}}$\\
    \hspace*{4ex}${\cLimSeqX {n} {\cSumX {i} {1} {n}
    {\cBin {\cFunAppX {f} {\cBin {a} {+} 
    {\cBin {\cFracX {\cBinX {b} {-} {a}} {n}} {*} {i}}}}
    {*} {\cFracX {\cBinX {b} {-} {a}} {n}}}}}$
    \hfill (definite integral).\\[2ex]

  $\mName{Thm27}$: ${\cForallBX {f,g} {\cFunTyX {R} {R}} {a,b} {R} {}}$\\
    \hspace*{2ex}${\cImpliesX {\cAnd 
    {\cFunAppCX {\textsf{cont-on-closed-int}} {f} {a} {b}}
    {\cEqX {g} {\cFunAbsX {x} {R} 
    {\cIntegralX {a} {x} {\cFunApp {f} {s}} {s}}}}} {}}$\\
    \hspace*{4ex}$({\cAndX
    {\cForall {x} {R} {\cImpliesX {\cBinBX {a} {<} {x} {<} {b}}
    {\cEqX {\cFunAppBX {\textsf{deriv-at}} {g} {x}} 
    {\cFunAppX {f} {x}}}}} {}}$\\
    \hspace*{5ex}${\cAndX 
    {\cEqX {\cFunAppBX {\textsf{right-deriv-at}} {g} {a}} 
    {\cFunAppX {f} {a}}} {}}$\\
    \hspace*{5ex}${\cEqX {\cFunAppBX {\textsf{left-deriv-at}} {g} {b}} 
    {\cFunAppX {f} {b}}})$
    \hfill (fundamental theorem of calculus).\\[2ex]

  $\mName{Thm28}$: ${\cForallBX {f,g} {\cFunTyX {R} {R}} {a,b} {R} {}}$\\
    \hspace*{2ex}${\cImpliesX {\cAnd
    {\cFunAppCX {\textsf{cont-on-closed-int}} {f} {a} {b}}
    {\cEqX {f} {\cFunAppX {\mName{deriv}} {g}}}}
    {\cEqX {\cIntegralX {a} {b} {\cFunApp {f} {x}} {x}}
    {\cBinX {\cFunAppX {g} {b}} {-} {\cFunAppX {g} {a}}}}}$\\
    \phantom{x} \hfill (corollary of fundamental theorem of calculus).\\[2ex]
\hline
\end{tabular}
\ec
\caption{Some Calculus Definitions and Theorems from a Development of the Alonzo Theory $\mathsf{COF}$}
\label{fig:COF-calculus}
\end{figure*}

\begin{figure*}
\bc
\begin{tabular}{|lll|}
\hline

  &&\phantom{x}\\

  $\cSum {\mathbf{i}} {\mathbf{M}_R} {\mathbf{N}_R} {\mathbf{A}_R}$
& stands for
& $\cFunAppCX {\mName{sum}_{\cFunTyCX {R} {R} {\cFunTy {R} {R}} {R}}} {\mathbf{M}_R}
  {\mathbf{N}_R} {\cFunAbs {\mathbf{i}} {R} {\mathbf{A}_R}}$.\\[3ex]

  $\cLim {\mathbf{x}} {\mathbf{A}_R} {\mathbf{B}_R}$
& stands for
& $\cFunAppBX {\mName{lim}_{\cFunTyBX {\cFunTy {R} {R}} {R} {R}}}
  {\cFunAbs {\mathbf{x}} {R} {\mathbf{B}_R}} {\mathbf{A}_R}$.\\[3ex]

  $\cLimSeq {\mathbf{n}} {\mathbf{B}_R}$
& stands for
& $\cFunAppX {\textsf{lim-seq}_{\cFunTyX {\cFunTy {R} {R}} {R}}}
  {\cFunAbs {\mathbf{n}} {N_{\cSetTy {R}}} {\mathbf{B}_R}}$.\\[3ex]

  $\cIntegral {\mathbf{A}_R} {\mathbf{B}_R} {\mathbf{C}_R} {\mathbf{x}}$
& stands for
& $\cFunAppCX {\mName{integral}_{\cFunTyCX {\cFunTy {R} {R}} {R} {R} {R}}}
  {\cFunAbs {\mathbf{x}} {R} {\mathbf{C}_R}} {\mathbf{A}_R} {\mathbf{B}_R}$.\\[3ex]
\hline
\end{tabular}
\ec
\caption{Some Notational Definitions for the Alonzo Theory $\mathsf{COF}$}
\label{fig:notational-defs}
\end{figure*}

There is a formal proof system for Alonzo~\cite{Farmer25} which is
derived from Andrews' elegant proof system for
{\qzero}~\cite{Andrews02}.  Proofs, however, are not required to be
formal in our implementation of the free approach, and so
requirement R2 is satisfied.

Alonzo has a module system for organizing mathematical knowledge using
the little theories method~\cite{Farmer25}.  (Since partial functions
naturally arise from theory morphisms~\cite{Farmer94}, the little
theories method works best with a logic like Alonzo that supports
partial functions.)  It includes modules for constructing theories and
theory morphisms, developing theories by defining new concepts and
posing and proving conjectures, and transporting definitions and
theorems from one theory to another.  Thus Alonzo fully satisfies
requirement R3.

Our implementation of the free approach currently provides only the
simplest level of software support: LaTeX macros for presenting Alonzo
types and expressions and LaTeX environments for presenting Alonzo
modules.  Other levels of software support are possible; see the
discussion in Chapter~16 of~\cite{Farmer25}.  Alonzo has not been
implemented in a proof assistant, but since it is closely related to
LUTINS, it could be implemented in much the same way that LUTINS is
implemented in IMPS.  Thus R4 is only partially satisfied now, but it
could be fully satisfied with the addition of more levels of software
support.

We have presented in~\cite{FarmerZvigelsky25} an example of a
substantial body of mathematical knowledge associated with monoid
theory that is formalized in Alonzo using the little theories method.
This example illustrates how mathematics can be done using the free
approach to formal mathematics.  The example is described in the next
section.

Our Alonzo-based implementation of the free approach demonstrates that
all the benefits of formal mathematics, except possibly mechanically
checked results, can be achieved when communication is the primary
focus and only minimal software support is utilized.  In other words,
it demonstrates that there is a path to formal mathematics that will
better serve most mathematics practitioners than the
prove-everything-with-a-proof-assistant path.

\section{An Example}\label{sec:example}

The example given in~\cite{FarmerZvigelsky25} mentioned above presents
a a body of mathematical knowledge about monoids as a theory graph in
Alonzo using the free approach.  A robust theory graph (see
Figure~\ref{fig:monoid-theory}) is constructed using Alonzo modules.
It contains 12 theories and 18 theory morphisms (8 inclusions and 10
noninclusions).  Each of the theories is developed by defining new
concepts and stating and proving theorems.  For example, the theorem
\textsf{id-elt-is-unique} mentioned above that says the identity
element of a monoid is unique is stated as the sentence
\[{\cForallX {x} {M} {\cImpliesX {\cForall {y} {M}
      {\cBinBX {\cBinX {x} {\cdot} {y}} {=} {\cBinX {y} {\cdot} {x}}
        {=} {y}}} {\cEqX {x} {\mathsf{e}}}}}\] of the language of the
theory $\mathsf{MON}$ presented using Alonzo's compact notation and
proved from the axioms of $\mathsf{MON}$.  The theory graph includes
15 definitions and 36 theorems.  The 51 proofs that validate the
definitions and the theorems are all included in the appendix
of~\cite{FarmerZvigelsky25} and are written in a traditional style,
but some make use of the axioms, rules of inference, and metatheorems
of Alonzo's formal proof system.  The formalization is done with only
the simplest level of software support, just the LaTeX macros and
environments mentioned~above.

The 12 theories in the theory graph play different roles.  For two
theories $T$ and $T'$ in the theory graph, we say \emph{$T$ provides
services to $T'$} and \emph{$T'$ receives services from $T$} when
definitions or theorems are transported from $T$ to $T'$ via a theory
morphism from $T$ to $T'$.  Two theories in the theory graph only
provide services, two only receive services, and four both provide and
receive services.

\textsf{FUN-COMP}, a theory of function composition, is an example of
the first kind of theory.  Theorems about general function composition
are transported from it to \textsf{ONE-BT}, a theory having just a
single base type.  The theories \textsf{FUN-COMP} and \textsf{ONE-BT}
have no constants nor axioms, yet they are still useful by virtue of
being theories in a richly endowed logic, i.e., Alonzo.

\textsf{ONE-BT-with-SC}, a theory of standard transformation monoids
that extends \textsf{ONE-BT}, is an example of the second kind of
theory.  Theorems are transported to it from \textsf{MON-ACT}, a
theory of monoid actions, that show a standard transformation monoid
exhibits the structure of monoid action.

\textsf{MON}, a theory of monoids that is central to the theory graph,
is an example of the third kind of theory.  Definitions and theorems
are transported from it to \textsf{ONE-BT} and \textsf{MON-HOM}, a
theory of monoid homomorphisms.  Theorems are transported to it from
\textsf{MON-HOM} and \textsf{MON-ACT}.  Theorems are also transported
from \textsf{MON} to itself via two different theory morphisms.

This example demonstrates four important things.  First, formal
mathematics can be performed with minimal software support and only
traditional proofs.  Second, the free approach to formal mathematics
facilitates communication and ease of use better than the standard
approach.  Third, the little theories method enables mathematical
knowledge to be formalized to maximize clarity and minimize
redundancy.  And fourth, Alonzo is well suited for expressing and
reasoning about mathematical ideas.

\section{Free Approach and Artificial Intelligence}

The reader may be wondering whether LLMs will, in time, make the free
approach obsolete.  Great strides have been made over the last few
years in using LLMs to perform \emph{autoformalization}, i.e., the
automatic translation of traditional mathematics to fully certified
formal mathematics.  For example, see~\cite{RammalEtAl26}.  It is a
relatively easy task to automatically generate a natural language
presentation of a formal proof that better communicates the proof.
This is just a matter of translating the usual notation for a formal
proof to a controlled language notation for the proof.  Thus it
appears that in the future it will be possible to automatically
translate traditional mathematics into formal mathematics presented
using natural language.  However, the resulting natural language
proofs will be just as much certification-oriented as the original
formal proofs.

What is really needed is \emph{autoexplanation}, i.e, the automatic
translation of certification-oriented formal mathematics to a
communication-oriented formal mathematics.  This is likely to be a
much bigger challenge than autoformalization since it requires a deep
understanding of the mathematical context in which the proof resides
and the target audience for the proof.  We thus expect that
autoformalization will be much more successful at certifying a
traditional proof than autoexplanation will be at communicating the
key ideas embodied in a formal proof.

The effective explanation of a proof for a human requires human
creativity and intelligence.  The effective presentation of a body of
mathematical knowledge for a human requires even more human creativity
and intelligence.  So, for the foreseeable future, it is not likely
that AI will not make the free approach~obsolete.

\section{Conclusion}

Formal mathematics, i.e., mathematics done within the framework of a
formal logic, offers huge benefits to the mathematics practitioner.
There are millions of mathematics practitioners --- and billions if
mathematics students are included --- yet only a minuscule portion of
these do mathematics with the help of a formal logic and supporting
software.  The free approach to formal mathematics has the potential
to enable a much larger portion of mathematics practitioners to reap
the benefits of formal mathematics than the standard approach in which
mathematics is done using a proof assistant and all details are
formally proved and mechanically checked.

The free approach can be an attractive option for mathematics
practitioners who have not employed formal mathematics due to the high
cost of the standard approach.  As we have illustrated
in~\cite{FarmerZvigelsky25}, the cost of employing the free approach
can be much lower than doing everything with a proof assistant.  The
free approach can also be an attractive option for mathematics
practitioners for whom communication is more important than
certification.  When the mathematics is well understood, as in
mathematics education and routine applications, there is much less
need to construct and mechanically check formal proofs.  Another
benefit of the free approach is that it can be a stepping stone to
help mathematics practitioners cross the void between traditional
mathematics and fully certified formal mathematics.

The great, and largely unrealized, promise of formal mathematics
cannot be achieved by pursuing just the standard approach to formal
mathematics.  The free approach is needed in addition.  Therefore, we
strongly encourage the mathematics community to develop logics,
software, and libraries of formal mathematical knowledge to support
the free approach and to train mathematics practitioners to use them.

We look forward to a time when mathematics practitioners routinely
express and reason about their ideas within a formal logic and high
school and university mathematics students routinely build bodies of
mathematical knowledge as theory graphs.

\bibliography{simple-type-theory-bib}
\bibliographystyle{plain}

\end{document}